\documentclass[10pt,fleqn,leqno]{article}
\usepackage{amsmath}

\title{Kahler-Einstein Structures of General Natural Lifted Type on the
Cotangent Bundles}
\author{S.~L.~Dru\c t\u a}

\date{}

\numberwithin{equation}{section}

\begin{document}
\maketitle

\begin{abstract} \normalsize\bf Abstract. \rm We study the conditions
under which the cotangent bundle $T^*M$ of a Riemaannian manifold
$(M,g)$, endowed with a K\"ahlerian structure $(G,J)$ of general
natural lift type (see \cite{Druta1}), is Einstein. We first
obtain a general natural K\"ahler-Einstein structure on the
cotangent bundle $T^*M$. In this case, a certain parameter,
$\lambda$ involved in the condition for $(T^*M,G,J)$ to be a
K\"ahlerian manifold, is expressed as a rational function of the
other two, the value of the constant sectional curvature, $c$, of
the base manifold $(M,g)$ and the constant $\rho$ involved in the
condition for the structure of being Einstein. This  expression of
$\lambda$ is just that involved in the condition for the
K\"ahlerian manifold to have constant holomorphic sectional
curvature (see \cite{Druta2}). In the second case, we obtain a
general natural K\"ahler-Einstein structure only on $T_0M$, the
bundle of nonzero cotangent vectors to $M$. For this structure,
$\lambda$ is expressed as another function of the other two
parameters, their derivatives, $c$  and $\rho$.
\end{abstract}
\vskip3mm {\bf Mathematics Subject Classification 2000:} primary
53C55, 53C15, 53C07.\\{\bf~ Key words:} cotangent bundle,
Riemannian metric, general natural lift, K\"ahler-Einstein
structure.

\normalsize \section{Introduction}
\renewcommand \theequation{\thesection.\arabic{equation}}
A few natural lifted structures introduced on the cotangent bundle
$T^*M$ of a Riemannian manifold $(M,g)$, have been studied in
recent papers such as \cite{Druta1}, \cite{Druta2},
\cite{Munteanu}, \cite{OprPap1}, \cite{OprPap2},
\cite{OprPap3}--\cite{Porosniuc2}. The similitude between some
results from the mentioned papers and results from the geometry of
the tangent bundle $TM$ (e.g. \cite{Oproiu4}, \cite{Oproiu3}), may
be explained by the duality cotangent bundle-tangent bundle. The
fundamental differences between the geometry of the cotangent
bundle and that of the tangent bundle of a Riemannian manifold,
are due to the  different construction of lifts to $T^*M$, which
cannot be defined just like in the case of $TM$ (see
\cite{YanoIsh}).

The results from \cite{Kolar}, \cite{KowalskiSek},
\cite{KrupkaJan}, \cite{Terng}, concerning the natural lifts, and
the classification of the natural vector fields on the tangent
bundle of a pseudo-Riemannian manifold, made by Jany$\check{s}$ka
in \cite{Janyska}, allowed the present author to introduce in the
paper \cite{Druta1}, a general natural almost complex structure
$J$ of lifted type on the cotangent bundle $T^*M$, and a general
natural lifted metric $G$ defined by the Riemannian metric $g$ on
$T^*M$ (see the paper \cite{Oproiu4} by Oproiu, for the case of
the tangent bundle). The main result from \cite{Druta1} is that
the family of general natural K\" ahler structures on $T^*M$
depends on three essential parameters (one is a certain
proportionality factor obtained from the condition for the
structure to be almost Hermitian and the other two are
coefficients involved in the definition of the integrable almost
complex structure $J$ on $T^*M$).

In the present paper we are interested in finding the conditions
under which the cotangent bundle $T^*M$ of a Riemannian manifold
$(M,g)$, endowed with a  K\"ahlerian structure $(G,J)$ of general
natural lift type (see \cite{Druta1}), is an Einstein manifold. To
this aim, we have to study the vanishing conditions for the
components of the difference between the Ricci tensor of
$(T^*M,G,J)$ and $\rho G$, where $\rho$ is a constant.

After some quite long computations with  the RICCI package from
the program Mathematica, we obtain two cases in which a general
natural K\"ahlerian manifold $(T^*M,G,J)$ is Einstein. In the
first case, $(T^*M,G,J)$ is a K\"ahler-Einstein manifold if the
proportionality factor $\lambda$, involved in the condition for
the manifold to be K\"ahlerian, is expressed as a rational
function of the first two essential parameters, their derivatives,
the values of the constant sectional curvature of the base
manifold $(M,g)$, and the constant $\rho$, from the condition for
the manifold to be Einstein. In this case the expression of
$\lambda$ leads to the condition obtained in \cite{Druta2} for
$(T^*M,G,J)$ to have constant holomorphic sectional curvature. In
the second case, $(G,J)$ is a K\"ahler-Einstein structure on the
the bundle of nonzero cotangent vectors to $M$, $T_0^*M$, if
$\lambda^\prime$ is expressed as a certain function of $\lambda$,
the other two parameters, their first order derivatives, $c$ and
$\rho$. The similar problem on tangent bundle $TM$ was treated by
Oproiu and Papaghiuc in the paper \cite{OprPap4}.

The manifolds, tensor fields and other geometric objects
considered in this paper are assumed to be differentiable of class
$C^\infty $ (i.e. smooth). The Einstein summation convention is
used throughout this paper, the range of the indices
$h,i,j,k,l,m,r $ being always $\{1,\dots ,n\}$.

\section{Preliminary results}

The cotangent bundle of a smooth $n$-dimensional Riemannian
manifold may be endowed with a structure of a $2n$-dimensional
smooth manifold, induced from the structure of the base manifold.
If $(M,g)$ is a smooth Riemannian manifold of the dimension $n$,
we denote its cotangent bundle by $\pi :T^*M\rightarrow M$. From
every local chart $(U,\varphi )=(U,x^1,\dots ,x^n)$  on $M$, it is
induced a local chart $(\pi^{-1}(U),\Phi )=(\pi^{-1}(U),q^1,\dots
, q^n,$ $p_1,\dots ,p_n)$, on $T^*M$, as follows. For a cotangent
vector $p\in \pi^{-1}(U)\subset T^*M$, the first $n$ local
coordinates $q^1,\dots ,q^n$ are  the local coordinates of its
base point $x=\pi (p)$ in the local chart $(U,\varphi )$ (in fact
we have $q^i=\pi ^* x^i=x^i\circ \pi, \ i=1,\dots n)$. The last
$n$ local coordinates $p_1,\dots ,p_n$ of $p\in \pi^{-1}(U)$ are
the vector space coordinates of $p$ with respect to the natural
basis $(dx^1_{\pi(p)},\dots , dx^n_{\pi(p)})$, defined by the
local chart $(U,\varphi )$,\ i.e. $p=p_idx^i_{\pi(p)}$.

The $M$-tensor fields on the cotangent bundle may be introduced in
the same manner as the $M$-tensor fields were introduced in the
paper \cite{Mok} on the tangent bundle of a Riemannian manifold.

On $T^*M$, a few useful $M$-tensor fields may be obtained as
follows. Let $v,w:[0,\infty ) \rightarrow {\bf R}$ be smooth
functions and let $\|p\|^2=g^{-1}_{\pi(p)}(p,p)$ be the square of
the norm of the cotangent vector $p\in \pi^{-1}(U)$ ($g^{-1}$ is
the tensor field of type (2,0) having the components $(g^{kl}(x))$
which are the entries of the inverse of the matrix $(g_{ij}(x))$
defined by the components of $g$ in the local chart $(U,\varphi
)$). The components $vg_{ij}(\pi(p))$, $p_i$, $w(\|p\|^2)p_ip_j $
define respective $M$-tensor fields of types $(0,2)$, $(0,1)$,
$(0,2)$ on $T^*M$. Similarly, the components $vg^{kl}(\pi(p))$,
$g^{0i}=p_hg^{hi}$, $w(\|p\|^2)g^{0k}g^{0l}$ define respective
$M$-tensor fields of type $(2,0)$, $(1,0)$, $(2,0)$ on $T^*M$. Of
course, all the components considered above are in the induced
local chart $(\pi^{-1}(U),\Phi)$.

We recall the splitting of the tangent bundle to $T^*M$ into the
vertical distribution $VT^*M= {\rm Ker}\ \pi _*$ and the
horizontal one determined by the Levi Civita connection $\dot
\nabla $ of $g$:
\begin{eqnarray}\label{descomp}
~~~~~~~~~~~~~~~~~~~~~~~~~~TT^*M=VT^*M\oplus HT^*M.
\end{eqnarray}
If $(\pi^{-1}(U),\Phi)=(\pi^{-1}(U),q^1,\dots ,q^n,p_1,\dots
,p_n)$ is a local chart on $T^*M$, induced from the local chart
$(U,\varphi )= (U,x^1,\dots ,x^n)$, the local vector fields
$\frac{\partial}{\partial p_1}, \dots , \frac{\partial}{\partial
p_n}$ on $\pi^{-1}(U)$ define a local frame for $VT^*M$ over $\pi
^{-1}(U)$ and the local vector fields $\frac{\delta}{\delta
q^1},\dots ,\frac{\delta}{\delta q^n}$ define a local frame for
$HT^*M$ over $\pi^{-1}(U)$, where
\begin{eqnarray*}
~~~~~~~~~~~~~~~~~~~~~~\frac{\delta}{\delta
q^i}=\frac{\partial}{\partial q^i}+\Gamma^0_{ih}
\frac{\partial}{\partial p_h},\ \ \ \Gamma ^0_{ih}=p_k\Gamma
^k_{ih},
\end{eqnarray*}
and $\Gamma ^k_{ih}(\pi(p))$ are the Christoffel symbols of $g$.

The set of vector fields $\{\frac{\partial}{\partial p_1},\dots
,\frac{\partial}{\partial p_n}, \frac{\delta}{\delta q^1},\dots
,\frac{\delta}{\delta q^n}\}$ defines a local frame on $T^*M$,
adapted to the direct sum decomposition (\ref{descomp}).

We consider
\begin{eqnarray*}
~~~~~~~~~~t=\frac{1}{2}\|p\|^2=\frac{1}{2}g^{-1}_{\pi(p)}(p,p)=\frac{1}{2}g^{ik}(x)p_ip_k,
\ \ \ p\in \pi^{-1}(U)
\end{eqnarray*}
the energy density defined by $g$ in the cotangent vector $p$. We
have $t\in [0,\infty)$ for all $p\in T^*M$.

The computations will be done in local coordinates, using a local
chart $(U,\varphi)$ on $M$ and the induced local chart
$(\pi^{-1}(U),\Phi)$ on $T^*M$.

We shall use the following lemma, which may be proved easily.

\bf{Lemma 2.1.} \it{If $n>1$ and $u,v$ are smooth functions on
$T^*M$ such that
\begin{eqnarray*}
u g_{ij}+v p_ip_j=0, \quad u g^{ij}+v g^{0i}g^{0j}=0,\quad or
\quad u\delta ^i_j+vg^{0i} p_j=0,\ \forall i,j=\overline{1,n},
\end{eqnarray*}
on the domain of any induced local chart on $T^*M$, then $u=0,\
v=0$.}

\rm In the paper \cite{Druta1}, the present author considered the
real valued smooth functions $a_1,\ a_2,\ a_3,$ $\ a_4,\ b_1,\
b_2,\ b_3,\ b_4$ on $[0,\infty)\subset {\bf R}$ and studied a
general natural tensor of type $(1,1)$ on $T^*M$, defined by the
relations
\begin{equation}\label{Jinvar}
\left\{
\begin{array}{l}
JX^H_p=a_1(t)
(g_X)^V_p+b_1(t)p(X)p_p^V+a_4(t)X_p^H+b_4(t)p(X)(p^\sharp)_p^H,
\\ \mbox{ }  \\
J\theta^V_p=a_3(t)\theta^V_p+b_3(t)g^{-1}_{\pi(p)}
(p,\theta)p_p^V-a_2(t)(\theta^\sharp)_p^H- b_2(t)g^{-1}_{\pi(p)}
(p,\theta)(p^\sharp)_p^H,
\end{array}
\right.
\end{equation}
in every point $p$ of the induced local card $(\pi^{-1}(U),\Phi)$
on $T^*M$, $\forall ~X \in \mathcal{X}(M), \forall~ \theta \in
\Lambda^1 (M)$, where $g_X$ is the 1-form on $M$ defined by
$g_X(Y)=g(X,Y),\ \forall Y\in \mathcal{X}(M)$,
$\theta^\sharp=g^{-1}_\theta$ is a vector field on $M$ defined by
$g(\theta^\sharp,Y)=\theta (Y),~\forall~ Y \in \mathcal{X}(M)$,
the vector $p^\sharp$ is tangent to $M$ in $\pi (p)$, $p^V$ is the
Liouville vector field on $T^*M$ , and $(p^\sharp)^H$ is the
similar horizontal vector field on $T^*M$.

The definition of the general natural lift given by
(\ref{Jinvar}), is based on the Jany$\check{s}$ka's classification
of the natural vector fields on the tangent bundle, but the
construction is different, being specific for the cotangent
bundle.

\bf{Theorem 2.1. \rm{(\cite{Druta1})\label{th4}} \it{A natural
tensor field $J$ of type $(1,1)$ on $T^*M$, given by
$(\ref{Jinvar})$, defines an almost complex structure on $T^*M$,
if and only if $a_4=-a_3, b_4=-b_3$ and the coefficients $a_1,\
a_2,\ a_3,\ b_1,\ b_2$ and $b_3$ are related by}
\begin{equation}\label{rel4}
~~~~~~~a_1a_2=1+a_3^2\ ,\ \ \
(a_1+2tb_1)(a_2+2tb_2)=1+(a_3+2tb_3)^2.
\end{equation}

\rm Studying the vanishing conditions for the Nijenhuis tensor
field $N_J$, we may state:

\bf{Theorem 2.2.} \rm(\cite{Druta1})\label{th3} \it{ Let $(M,g)$
be an $n(>2)$-dimensional connected Riemannian manifold. The
almost complex structure $J$ defined by {\rm(\ref{Jinvar})} on
$T^*M$ is integrable if and only if $(M,g)$ has constant sectional
curvature $c$ and the  coefficients $b_1,\ b_2,\ b_3$ are given
by:}
\begin{equation}
~~~~~~~~~~~~~~~~~~~\left\{\begin{array}{ll}\label{integrab}
b_1=\frac{2 c^2 t a_2^2+2 c t a_1 a_2^\prime+a_1 a_1^\prime -c+3 c
a_3^2}{a_1-2 t a_1^\prime-2 c t a_2-4 c t^2
a_2^\prime},\\
b_2=\frac{2 t a_3^{\prime 2}- 2 t a_1^\prime a_2^\prime+c a_2^2+2
c t a_2 a_2^\prime+a_1 a_2^\prime}{a_1-2 t a_1^\prime-2 c t a_2-4
c t^2 a_2^\prime},\\
b_3 =\frac{a_1 a_ 3^\prime+ 2 c a_2 a_3+4 c t a_2^\prime a_3-
 2 c t a_ 2 a_3^\prime}{a_1-2 t a_1^\prime-2 c t a_2-4 c t^2 a_2^\prime}.\\
\end{array}
\right.
\end{equation}

\bf Remark 2.3. \rm The integrability conditions (\ref{integrab})
for the almost complex structure $J$ on $T^*M$, may be expressed
in the equivalent form
\begin{eqnarray*}
 ~~~~~~~~~~~~~~~~~~~~~ \left\{%
\begin{array}{ll}\label{integr}
  a_1^{\prime}=\frac{1}{a_1+2tb_1}(a_1b_1+c-3ca_3^2-4cta_3b_3), \\
    a_2^{\prime}=\frac{1}{a_1+2tb_1}(2a_3b_3-a_2b_1-ca_2^2), \\
    a_3^{\prime}=\frac{1}{a_1+2tb_1}(a_1b_3-2ca_2a_3-2cta_2b_3). \\
\end{array}%
\right.
\end{eqnarray*}

\rm In the paper \cite{Druta1}, the author defined a Riemannian
metric $G$ of general natural lift type, given by the relations
\begin{equation}\label{Ginvar}
\left\{
\begin{array}{l}
G_p(X^H, Y^H) = c_1(t)g_{\pi(p)}(X,Y) + d_1(t)p(X)p(Y),
\\ \mbox{ }\\
G_p(\theta^V,\omega^V) = c_2(t)g^{-1}_{\pi(p)}(\theta,\omega) +
d_2(t)g^{-1}_{\pi(p)}(p,\theta)g^{-1}_{\pi(p)}(p,\omega),
\\ \mbox{ }\\
G_p(X^H,\theta^V) = G_p(\theta^V,X^H)
=c_3(t)\theta(X)+d_3(t)p(X)g^{-1}_{\pi(p)}(p,\theta),
\end{array}
\right.
\end{equation}
$\forall~ X,Y \in \mathcal{X}(M),$ $\forall~ \theta, \omega \in
\Lambda^1(M), \forall~p \in T^*M$.

\vskip3mm

Using the adapted frame $\{\frac{\partial}{\partial
p_i},\frac{\delta}{\delta q^j}\}_{i,j=1,\dots,n} $ on $T^*M$, we
may write the expression (\ref{Ginvar}) in the next form
\begin{equation}\label{rel11}
\left\{
\begin{array}{l}
G(\frac{\delta}{\delta q^i}, \frac{\delta}{\delta
q^j})=c_1(t)g_{ij}+ d_1(t)p_ip_j=G^{(1)}_{ij},
\\   \mbox{ } \\
G(\frac{\partial}{\partial p_i}, \frac{\partial}{\partial
p_j})=c_2(t)g^{ij}+d_2(t)g^{0i}g^{0j}=G_{(2)}^{ij},
\\   \mbox{ } \\
G(\frac{\partial}{\partial p_i},\frac{\delta}{\delta q^j})=
G(\frac{\delta}{\delta q^i},\frac{\partial}{\partial p_j})=
c_3(t)\delta_i^j+d_3(t)p_ig^{0j}=G3_i^j,
\end{array}
\right.
\end{equation}
where $c_1,\ c_2,\ c_3,\ d_1,\ d_2,\ d_3$ are six smooth functions
of the density energy on $T^*M$. The conditions for $G$ to be
positive definite are assured if
\begin{eqnarray*}
c_1+2td_1>0,\quad c_2+2td_2>0,\quad
(c_1+2td_1)(c_2+2td_2)-(c_3+2td_3)^2>0.
\end{eqnarray*}

The author proved the following result

\bf{Theorem 2.3.} \rm(\cite{Druta1})\it{\label{th4} {The family of
Riemannian metrics  $G$ of general natural lifted type on $T^*M$
such that $(T^*M,G,J)$ is an almost Hermitian manifold, is given
by {\rm(\ref{rel11})}, provided that the coefficients $c_1,\ c_2,\
c_3,\  d_1,\ d_2,$ and $d_3$ are related to the coefficients
$a_1,\ a_2,\ a_3,\ b_1,\ b_2,$ and $b_3$ by the following
proportionality relations
\begin{eqnarray*}
~~~~\frac{c_1}{a_1} =\frac{c_2}{a_2}=\frac{c_3}{a_3} =
\lambda,\qquad \frac{c_1+2 t d_1}{a_1+2 t b_1} =\frac{c_2+2 t
d_2}{a_2+2 t b_2}=\frac{c_3+2 t d_3}{a_3+2 t b_3} = \lambda +2 t
\mu,
\end{eqnarray*}
where the proportionality coefficients $\lambda>0 $ and $\lambda
+2t\mu>0$ are functions depending on $t$.}

\rm Considering the two-form $\Omega $ defined by the almost
Hermitian structure $(G,J)$ on $T^*M$, given by $\Omega
(X,Y)=G(X,JY),$ for any vector fields $X,Y$ on $T^*M$, we may
formulate the main results from \cite{Druta1}:

\bf{Theorem 2.4.} \rm(\cite{Druta1}) \it{\label{th6} The almost
Hermitian structure $(T^*M,G,J)$ is almost K\"{a}hlerian if and
only if}
$$
\mu=\lambda ^\prime .
$$

\bf{Theorem 2.5.} \it{A general natural lifted almost Hermitian
structure $(G,J)$ on $T^*M$ is K\"ahlerian if and only if  the
almost complex structure $J$ is integrable (see Theorem \rm 2.2)
\it and $\mu=\lambda^\prime$.}

\rm Examples of such structures may be found in \cite{OprPor1},
\cite{Porosniuc3}.

\section{General natural K\"ahler-Einstein structures on  cotangent bundles}

The Levi-Civita connection $\nabla$ of the Riemannian manifold
$(T^*M,G)$ is obtained from the Koszul formula, and it is
characterized by the conditions
$$
\nabla G=0,\ T=0,
$$
where $T$ is the torsion tensor of $\nabla.$

\vskip2mm

In the case of the cotangent bundle $T^*M$ we may obtain the
explicit expression of $\nabla$.

The symmetric $2n\times 2n$ matrix
$$
\left(%
\begin{array}{cc}
 G^{(1)}_{ij} & G3^j_i \\
  G3^i_j & G_{(2)}^{ij} \\
\end{array}%
\right)
$$
associated to the metric $G$ in the base $\{\frac{\delta}{\delta
q^i},\frac{\partial}{\partial p_j}\}_{i,j=1,\dots,n}$ has the
inverse
$$
\left(%
\begin{array}{cc}
  H_{(1)}^{ij} & H3_i^j \\
  H3^i_j & H^{(2)}_{ij} \\
\end{array}%
\right)
$$
where the entries are the blocks
\begin{equation}\label{matrinv}
H_{(1)}^{kl}=e_1g^{kl}+f_1g^{0k}g^{0l},\quad
H^{(2)}_{kl}=e_2g_{kl}+f_2p_kp_l,\quad
H3^k_l=e_3\delta^k_l+f_3g^{0k}p_l.
\end{equation}
Here $g^{kl}$ are the components of the inverse of the matrix
$(g_{ij})$, $g^{0k}=p_ig^{ik}$, and $e_1,\ f_1,\ e_2,\ f_2,$ $\
e_3$, $f_3:[0,\infty)\rightarrow \mathbf{R},$ some real smooth
functions. Their expressions are obtained  by solving the system:
$$
\left\{%
\begin{array}{ll}\label{sistinv}
   G^{(1)}_{ih}H_{(1)}^{hk}+G3_i^hH3^k_h=\delta_i^k,\\
    G^{(1)}_{ih}H3^h_k+G3_i^hH^{(2)}_{hk}=0,\\
    G3^i_hH_{(1)}^{hk}+G_{(2)}^{ih}H3_h^k=0,\\
    G3^i_hH3^h_k+G_{(2)}^{ih}H^{(2)}_{hk}=\delta^i_k,\\
\end{array}%
\right.
$$
in which we substitute the relations (\ref{rel11}) and
(\ref{matrinv}). By using Lemma 2.1,  we get $e_1,\ e_2,\ e_3$ as
functions of $c_1,\ c_2,\ c_3$
\begin{eqnarray}\label{inversa1}
~~~~~~~~e_1=\frac{c_2}{c_1c_2-c_3^2},\ \
e_2=\frac{c_1}{c_1c_2-c_3^2},\ \ e_3=-\frac{c_3}{c_1c_2-c_3^2},
\end{eqnarray}
and $f_1,\ f_2,\ f_3$ as functions of $c_1,\ c_2,\ c_3,$ $d_1,\
d_2,\ d_3,$ $e_1,\ e_2,\ e_3$

\begin{equation}\label{inversa2}
\begin{array}{c}
  f_1=-\frac{c_2d_1e_1 - c_3d_3e_1 - c_3d_2e_3 + c_2d_3e_3 +
2d_1d_2e_1t - 2d_3^2e_1t}{c_1c_2 - c_3^2 + 2c_2d_1t + 2c_1d_2t -
4c_3d_3t + 4d_1d_2t^2 - 4d_3^2t^2},\\\\
 f_2=\frac{(c_3 + 2d_3t)[(d_3e_1 + d_2e_3)(c_1 + 2d_1t) - (d_1e_1 +
d_3e_3)(c_3 + 2d_3t)]}{(c_2 + 2d_2t)[(c_1 + 2d_1t)(c_2 + 2d_2t) -
(c_3 + 2d_3t)^2]}-\frac{d_2e_2 + d_3e_3}{c_2 + 2d_2t},\\\\
 f_3=-\frac{(d_3e_1 + d_2e_3)(c_1 + 2d_1t) - (d_1e_1 + d_3e_3)(c_3
+ 2d_3t)}{(c_1 + 2d_1t)(c_2 + 2d_2t) - (c_3 + 2d_3t)^2}.
\end{array}
\end{equation}

Next we may obtain the expression of the Levi Civita connection of
the Riemannian metric $G$ on $T^*M$.

\vskip2mm

\bf{Theorem 3.1.} \it{The Levi-Civita connection $\nabla$ of\ $G$
has the following expression in the local adapted frame
$\{\frac{\delta}{\delta q^i},\frac{\partial}{\partial
p_j}\}_{i,j=1,\dots,n}$
$$
\left\{%
\begin{array}{ll}
   \displaystyle\nabla_{\frac{\partial}{\partial p_i}}
\frac{\partial}{\partial p_j}=Q^{ij}_{\ \
h}\frac{\partial}{\partial
p_h}+\widetilde{Q}^{ijh}\frac{\delta}{\delta q^h}, \qquad
\nabla_{\frac{\delta}{\delta q^i}} \frac{\partial}{\partial
p_j}=(-\Gamma^j_{ih}+\widetilde{P}_{i\ \ h}^{\
j})\frac{\partial}{\partial p_h}+P_i^{\ jh}\frac{\delta}{\delta
q^h}, \\
    \displaystyle \nabla_{\frac{\partial}{\partial p_i}}
\frac{\delta}{\delta q^j}=P_j^{\ ih}\frac{\delta}{\delta
q^h}+\widetilde{P}_{j\ h}^{\ i}\frac{\partial}{\partial p_h},
\qquad
    \nabla_{\frac{\delta}{\delta q^i}} \frac{\delta}{\delta
q^j}=(\Gamma^h_{ij}+\widetilde{S}_{ij}^{\ \
h})\frac{\delta}{\delta p_h}+S_{ijh}\frac{\partial}{\partial p_h},
\end{array}%
\right.
$$
where $\Gamma^h_{ij}$ are the Christoffel symbols of the
connection $\dot\nabla$ and $M$-tensor fields appearing as
coefficients in the above expressions are given as
\begin{equation}\label{PQS}
\left\{%
\begin{array}{ll}
    Q^{ij}_{\ \
h}=\frac{1}{2}(\partial^iG^{jk}_{(2)}+\partial^jG^{ik}_{(2)}-
\partial^kG^{ij}_{(2)})H^{(2)}_{kh}+\frac{1}{2}(\partial^iG3^j_k+
\partial^j G3^i_k) H3^k_h,\\
   \widetilde{Q}^{ijh}=\frac{1}{2}(\partial^iG^{jk}_{(2)}+\partial^jG^{ik}_{(2)}-
\partial^kG^{ij}_{(2)})H3_k^h+\frac{1}{2}(\partial^iG3^j_k+
\partial^jG3^i_k)H_{(1)}^{kh},\\\\
    P_j^{\ ih}=\frac{1}{2}(\partial^iG3^k_j-
\partial^kG3^i_j)H3^h_k+\frac{1}{2}(\partial^iG_{jk}^{(1)}-
R^0_{ljk}G_{(2)}^{li})H_{(1)}^{kh}, \\
    \widetilde{P}_{j\ h}^{\ i}=\frac{1}{2}(\partial^iG3^k_j-
\partial^kG3^i_j)H^{(2)}_{kh}+\frac{1}{2}(\partial^iG_{jk}^{(1)}-
R^0_{ljk}G_{(2)}^{li})H3^k_h, \\\\
   S_{ijh}=\frac{1}{2}(c_2R^{0}_{lij}-\partial^kG_{ij}^{(1)})H^{(2)}_{kh}-c_3R^0_{ijk}H3^k_h, \\
    \widetilde{S}^{\ \
h}_{ij}=\frac{1}{2}(c_2R^{0}_{lij}-\partial^kG_{ij}^{(1)})H3^h_k-c_3R^0_{ijk}H_{(1)}^{kh},  \\
\end{array}%
\right.
\end{equation}
where $R^h_{kij}$ are the components of the curvature tensor field
of the Levi Civita connection $\dot \nabla$ of the base manifold
$(M,g)$.} \vskip2mm

 \rm If we replace in (\ref{PQS}) the relations
(\ref{rel11}), which define the metric $G$, the expressions
(\ref{matrinv}) for the inverse matrix $H$ of $G$, and the
formulas (\ref{inversa1}), (\ref{inversa2}) we obtain the detailed
expressions of $P_i^{\ jh},\ Q^{ij}_{\ \ h},\ S_{ijh},\
\widetilde{P}_{j\ h}^{\ i},\ \widetilde{Q}^{ijh},\
\widetilde{S}^{\ \ h}_{ij}.$

\vskip2mm

The curvature tensor field $K$ of the connection $\nabla$ is
defined by
$$
K(X,Y)Z=\nabla_X\nabla_YZ-\nabla_Y\nabla_XZ-\nabla_{[X,Y]}Z,\ \ \
X,Y,Z\in \Gamma (TM).
$$

By using the local adapted frame $\{\frac{\delta}{\delta
q^i},\frac{\partial}{\partial p_j}\}_{i,j=1,\dots
,n}=\{\delta_i,\partial^j\}_{i,j=1,\dots ,n}$ we obtain the
horizontal and vertical components of the curvature tensor field:
$$
K(\delta_i,\delta_j)\delta_k={QQQQ_{ijk}}^h\delta_h+QQQP_{ijkh}\partial^h,
$$
$$
K(\delta_i,\delta_j)\partial^k={QQPQ_{ij}}^{kh}\delta_h+QQPP_{ij\
\ h}^{\ \ k}\partial^h,
$$

$$
K(\partial^i,\partial^j)\delta_k=PPQQ^{ij\ \ h }_{\ \
k}\delta_h+{PPQP^{ij}}_{kh}\partial^h,
$$
$$
K(\partial^i,\partial^j)\partial^k=PPPQ^{ijkh}\delta_h+{PPPP^{ijk}}_h\partial^h,
$$

$$
K(\partial^i,\delta_j)\delta_k=PQQQ^{i\ \ \ h}_{\
jk}\delta_h+{PQQP^i}_{jkh}\partial^h,
$$
$$
K(\partial^i,\delta_j)\partial^k=PQPQ^{i\ \ kh}_{\
j}\delta_h+PQPP^{i\ \ k}_{\ j\ \ h}\partial^h,
$$
where the coefficients are the $M$-tensor fields given by
$$
{QQQQ_{ijk}}^h=\widetilde{S}^{\ \ l}_{jk}\widetilde{S}^{\ \
h}_{il}+P^{\ lh}_iS_{jkl} -\widetilde{S}^{\ \
h}_{jl}\widetilde{S}^{\ \ l}_{ik}-P^{\ lh}_jS_{ikl}-R_{lij}^0P^{\
lh}_k+R^h_{kij},
$$
$$
QQQP_{ijkh}=\widetilde{S}^{\ \ l}_{jk}S_{ilh}+\widetilde{P}^ {\
l}_{i\ h}S_{jkl}- \widetilde{S}^{\ \
l}_{ik}S_{jlh}-\widetilde{P}^{\ l}_{j\ h}S_{ikl}-\widetilde{P}^{\
l}_{k\ h}R^0_{lij},
$$

$$
{QQPQ_{ij}}^{kh}=\widetilde{P}^{\ k}_{j\ \ l}P^{\ lh}_i+P_j^{\
kl}\widetilde{S}^{\ \ h}_{il}- \widetilde{P}^{\ k}_{i\ \ l}P_j^{\
lh}-P^{\ kl}_{i}\widetilde{S}^{\ \
h}_{jl}-R^0_{lij}\widetilde{Q}^{lkh},
$$
$$
QQPP_{ij\ \ h}^{\ \ k}\partial^h=\widetilde{P}^{\ k}_{j\ \
l}\widetilde{P}^{\ l}_{i\ \ h}+P_j^{\ kl}S_{ilh}- \widetilde{P}^{\
k}_{i\ \ l}\widetilde{P}^{\ l}_{j\ h}-P^{\
kl}_iS_{jlh}-R^0_{lij}Q_{\ \ h}^{lk}-R^k_{lij},
$$

$$
PPQQ^{ij\ \ h }_{\ \ k}\delta_h=\partial^iP_k^{\
jh}-\partial^jP_k^{\ ih}+\widetilde{P}^ {\ j}_{k\ \
l}\widetilde{Q}^{ilh}+ P_k^{\ jl}P_l^{\ ih}-\widetilde{P}^{\ i}_{
k\ l}\widetilde{Q}^{jlh}-P_k^{\ il}P_l^{\ jh},
$$
$$
{PPQP^{ij}}_{kh}=\partial^i\widetilde{P}^{\ j}_{ k\
h}-\partial^j\widetilde{P}^{\ i}_{k\ h}+\widetilde{P}^{\ j}_{k\
l}Q_{\ \ h}^{il} +P_k^{\ jl}\widetilde{P}^{\ i}_{ l\
h}-\widetilde{P}^{\ i}_{k\ l}Q_{\ \ h}^{jl}-P_k^{\
il}\widetilde{P}^{\ j}_{l\ \ h},
$$

$$
PPPQ^{ijkh}=\partial^i\widetilde{Q}^{jkh}-\partial^j\widetilde{Q}^{ikh}+
Q_{\ \ l}^{jk}\widetilde{Q}^{ilh}+\widetilde{Q}^{jkl}P^{\
ih}_l-Q_{\ \ l}^{ik}\widetilde{Q}^{jlh} -\widetilde{Q}^{ikl}P_l^{\
jh},
$$
$$
{PPPP^{ijk}}_h=\partial^iQ_{\ \ h}^{jk}-\partial^jQ_{\ \
h}^{ik}+Q_{\ \ l}^{jk}Q_{\ \
h}^{il}+\widetilde{Q}^{jkl}\widetilde{P}^{\ i}_{l\ \ h}- Q^{ik}_{\
\ l}Q_{\ \ h}^{jl}-\widetilde{Q}^{ikl}\widetilde{P}^{\ j}_{l\ \
h},
$$

$$
PQQQ^{i\ \ \ h}_{\ jk}\delta_h=\partial^i\widetilde{S}^{\ \
h}_{jk}+S_{jkl}\widetilde{Q}^{ilh}+\widetilde{S}^{\ \
l}_{jk}P_l^{\ ih} -\widetilde{P}^{\ i}_{k\ l}P_j^{\ lh}-P^{\
il}_k\widetilde{S}^{\ \ h}_{jl},
$$
$$
{PQQP^i}_{jkh}=\partial^iS_{jkh}+\widetilde{S}^{\ \ l}_{jk}Q_{\ \
h}^{il}+\widetilde{S}^{\ \ l}_{jk}\widetilde{P}^{\ i}_{l\ h}
-\widetilde{P}_{k\ l}^{\ i}\widetilde{P}^{\ l}_{j\ h}-P_k^{\
il}S_{jlh},
$$

$$
PQPQ^{i\ \ kh}_{\ j}=\partial^iP_j^{\ kh}+\widetilde{P}^{\ k}_{j\
l}\widetilde{Q}^{ilh}+P_j^{\ kl}P_l^{\ ih}- Q_{\ \ l}^{ik}P_j^{\
lh}-\widetilde{Q}^{ikl}\widetilde{S}^{\ \ h}_{jl},
$$
$$
PQPP^{i\ \ k}_{\ j\ \ \ h}=\partial^i\widetilde{P}^{\ k}_{j\
h}+\widetilde{P}^{\ k}_{j\ l}Q_{\ \ h}^{il}+P_j^{\
kl}\widetilde{P}^{\ i}_{l\ h}- Q_{\ \ l}^{ik}\widetilde{P}^{\
l}_{j\ h}-\widetilde{Q}^{ikl}S_{jlh}.
$$

In order to get the final expressions of the above $M$-tensor
fields, we have to compute the first and second order partial
derivatives with respect to the cotangential coordinates, $p_i$ of
the usual tensor fields involved in the definition of the
Riemannian metric $G$.

$$
\partial^iG^{(1)}_{jk}=c_1'g^{0i}g_{jk}+d_1'g^{0i}p_jp_k+
d_1 \delta^i_jp_k+d_1 p_j\delta^i_k,
$$
$$
\partial^iG_{(2)}^{jk}=c_2'g^{0i}g^{jk}+d_2'g^{0i}g^{0j}g^{0k}+
d_2 g^{ij}g^{0k}+d_2 g^{0j}g^{ik},
$$
$$
\partial^iG3^j_k=c_3'g^{0i}\delta^j_k+d_3'g^{0i}g^{0j}p_k+
d_3 g^{ij}p_k+d_3 g^{0j}\delta^i_k,
$$
$$
\partial^i\partial^jG^{(1)}_{kl}=c_1''g^{0i}g^{0j}g_{kl}+c_1'g^{ij}g_{kl}+
d_1''g^{0i}g^{0j}p_kp_l+ d_1'g^{ij}p_kp_l+
$$
$$
+d_1'g^{0j}\delta^i_kp_l+d_1'g^{0j}p_k\delta^i_l+d_1'g^{0i}\delta^j_kp_l+
d_1 \delta^j_k\delta^i_l+
$$
$$
+d_1'g^{0i}p_k\delta^j_l+d_1 \delta^i_k\delta^j_l,
$$
$$
\partial^i\partial^jG_{(2)}^{kl}=c_2''g^{0i}g^{0j}g^{kl}+c_2'g^{ij}g^{kl}+
d_2''g^{0i}g^{0j}g^{0k}g^{0l}+ d_2'g^{ij}g^{0k}g^{0l}+
$$
$$
+d_2'g^{0j}g^{ik}g^{0l}+d_2'g^{0j}g^{0k}g^{il}+d_2'g^{0i}g^{jk}g^{0l}+
d_2 g^{jk}g^{il}+
$$
$$
+d_2'g^{0i}g^{0k}g^{jl}+d_2 g^{ik}g^{jl},
$$
$$
\partial^i\partial^jG3^k_l=c_3''g^{0i}g^{0j}\delta^k_l+c_3'g^{ij}\delta^k_l+
d_3''g^{0i}g^{0j}g^{0k}p_l+ d_3'g^{ij}g^{0k}p_l+
$$
$$
+d_3'g^{0j}g^{ik}p_l+d_3'g^{0j}g^{0k}\delta^i_l+d_3'g^{0i}g^{jk}p_l+
d_3 g^{jk}\delta^i_l+
$$
$$
+d_3'g^{0i}g^{0k}\delta^j_l+d_3 g^{ik}\delta^j_l,
$$
$$
\partial^iH_{(1)}^{jk}=e_1'g^{0i}g^{jk}+f_1'g^{0i}g^{0j}g^{0k}+
f_1g^{ij}g^{0k}+f_1 g^{0j}g^{ik},
$$
$$
\partial^iH^{(2)}_{jk}=e_2'g^{0i}g_{jk}+f_2'g^{0i}p_jp_k+
f_2\delta^i_jp_k+f_2 p_j\delta^i_k,
$$
$$
\partial^iH3^j_k=e_3'g^{0i}\delta^j_k+f_3'g^{0i}g^{0j}p_k+
f_3g^{ij}p_k+f_3 g^{0j}\delta^i_k.
$$

We get the first order partial derivatives of the $M$-tensor
fields $P_i^{\ jh},\ Q^{ij}_{\ \ h},\ S_{ijh},$ $\widetilde{P}_{j\
h}^{\ i},$ $\widetilde{Q}^{ijh},$ $\widetilde{S}^{\ \ h}_{ij}$
with respect to the cotangential coordinates $p_i$  and we replace
these derivatives, and the expressions (\ref{inversa1}),
(\ref{inversa2}) of the functions $e_1,\ e_2,\ e_3,\ f_1,\ f_2,\
f_3$ and of their derivatives in order to obtain the components of
the curvature tensor as functions of $a_1,\ a_2,\ a_3$ and their
derivatives  of first, second and  third order only. The
expressions are obtained by using the Mathematica package RICCI.
$$
\partial^iQ_{\ \ h}^{jk}=\frac{1}{2}\partial^iH^{(2)}_{lh}(\partial^jG^{(2)}_{kl}+
\frac{1}{2}H^{(2)}_{lh}(\partial^i\partial^jG^{(2)}_{kl}+\partial^i\partial^kG_{(2)}^{jl}
-\partial^i\partial^lG_{(2)}^{jk})+
$$
$$
+\frac{1}{2}\partial^iH3^l_h(\partial^jG3^k_l+\partial^k G3^k_l)+
\frac{1}{2}H3^l_h(\partial^i\partial^jG3^k_l+\partial^i\partial^k
G3^k_l),
$$
$$
\partial_i\widetilde{Q}^{jkh}=\frac{1}{2}\partial^iH3_l^h(\partial^jG^{kl}_{(2)}+\partial^kG^{jl}_{(2)}-
\partial^lG^{jk}_{(2)})+
$$
$$
+\frac{1}{2}H3_l^h(\partial^i\partial^jG^{kl}_{(2)}+\partial^i\partial^kG^{jl}_{(2)}
-\partial^i\partial^lG^{jk}_{(2)})+
$$
$$
+\frac{1}{2}\partial^iH_{(1)}^{lh}(\partial^jG3^k_l+
\partial^kG3^j_l)+\frac{1}{2}H_{(1)}^{lh}(\partial^i\partial^jG3^k_l+\partial^i\partial^kG3^j_l),
$$
$$
\partial^i\widetilde{P}_{j\ h}^{\ k}=\frac{1}{2}\partial^iH^{(2)}_{lh}(\partial^kG3^l_j-
\partial^lG3^k_j)+ \frac{1}{2}H^{(2)}_{lh}(\partial^i\partial^kG3^l_j
-\partial^i\partial^lG3^k_j)+
$$
$$
+\frac{1}{2}\partial^iH3^l_h(\partial^kG_{jl}^{(1)}-
R^0_{mjl}G_{(2)}^{mk})+
$$
$$
+\frac{1}{2}H3^l_h(\partial^i\partial^kG_{jl}^{(1)}-
R^i_{mjl}G_{(2)}^{mk}-R^0_{mjl}\partial^iG_{(2)}^{mk}),
$$
$$
\partial^iP_j^{\ kh}=\frac{1}{2}\partial^iH3^h_l(\partial^kG3^l_j-
\partial^lG3^k_j)+
\frac{1}{2}H3^h_l(\partial^i\partial^kG3^l_j
-\partial^i\partial^lG3^k_j)+
$$
$$
+\frac{1}{2}\partial^iH_{(1)}^{hl}(\partial^kG_{jl}^{(1)}-
R^0_{mjl}G_{(2)}^{mk})+
$$
$$
+\frac{1}{2}H_{(1)}^{hl}(\partial^i\partial^kG_{jl}^{(1)}-
R^i_{mjl}G_{(2)}^{mk}-R^0_{mjl}\partial^iG_{(2)}^{mk}),
$$
$$
\partial^iS_{jkh}=\frac{1}{2}[(c_2^{\prime}p^iR^{0}_{mjk}+c_2R^i_{mjk}-\partial^i\partial^lG^{(1)}_{jk})H^{(2)}_{lh}+
(c_2R^{0}_{mjk}-\partial^lG_{jk}^{(1)})\partial^iH^{(2)}_{lh}]-
$$
$$
-c_3^{\prime}p^iR^0_{jkl}H3^l_h-c_3(R^i_{jkl}H3^l_h+R^0_{jkl}\partial^iH3^l_h),
$$
$$
\partial^i\widetilde{S}^{\ \ h}_{jk}=\frac{1}{2}[(c_2^{\prime}p^iR^{0}_{mjk}+c_2R^i_{mjk}-\partial^i\partial^lG^{(1)}_{jk})H3^h_l+
(c_2R^{0}_{mjk}-\partial^lG_{jk}^{(1)})\partial^iH3^h_l]-
$$
$$
-c_3^{\prime}p^iR^0_{jkl}H_{(1)}^{lh}-c_3(R^i_{jkl}H_{(1)}^{lh}+R^0_{jkl}\partial^iH_{(1)}^{lh}).
$$

\vskip3mm

In the following, we shall obtain the conditions under which the
general natural K\"ahlerian manifold $(T^*M,G,J)$ is an Einstein
manifold. The components of the Ricci tensor
$Ric(Y,Z)=trace(X\rightarrow K(X,Y)Z)$ of the K\"alerian manifold
$(T^*M,G,J)$ are given by the formulas:
$$
RicQQ_{jk} = Ric(\delta_j,\delta_k)= QQQQ^{\ \ \ \ h}_{hjk}
+PQQP^h_{\ jkh},
$$
$$
RicPP^{jk} = Ric(\partial_j,\partial_k)= PPPP^{hjk}_{\ \ \ \
h}-PQPQ^{j\ kh}_{\ h},
$$
$$
RicQP_j^{\ k}=Ric(\delta_j,\partial^k)=RicPQ^k_{\ j}=
Ric(\partial^k,\delta_j) = PQPP_{\ j \ \ h}^{h\ k}+QQPQ_{hj}^{\ \
\ kh}.
$$

The conditions for the general natural K\"ahlerian manifold
$(T^*M,G,J)$ to be Einstein, are
$$
\left\{%
\begin{array}{ll}
   RicQQ_{jk}-\rho G^{(1)}_{jk}=0,\\
   RicPP^{jk}-\rho G_{(2)}^{jk}=0,\\
   RicQP_j^{\ k}-\rho G3^k_j=0,\\
\end{array}%
\right.
$$
where $\rho$ is a constant.

\vskip2mm

After a straightforward computation, using the RICCI package from
Mathematica, the three differences which we have to study, become
of the next forms:
$$
\left\{%
\begin{array}{ll}
   RicQQ_{jk}-\rho G^{(1)}_{jk}=(\lambda + 2\lambda't)^2[\lambda(\lambda + 2\lambda't)\alpha_1 g_{jk} + \beta_1 p_j p_k],\\
   RicPP^{jk}-\rho G_{(2)}^{jk}=\lambda(\lambda + 2\lambda't)^2[(\lambda + 2\lambda't)\alpha_2 g^{jk} + 2\lambda\beta_2 g^{0j} g^{0k}],\\
   RicQP_j^{\ k}-\rho G3^k_j=(\lambda + 2\lambda't)^2[\lambda(\lambda
+ 2\lambda't)\alpha_3 \delta_j^k + \beta_3 p_j g^{0k}],\\
\end{array}%
\right.
$$
where $\alpha_1,\ \alpha_2,\ \alpha_3,\ \beta_1,\ \beta_2, \
\beta_3$ are rational functions depending on $a_1,\ a_3,\
\lambda,$ their derivatives of the first two orders, and $\rho$.
We do not present here the explicit expressions of the functions,
since they are quite long.

\vskip2mm

Using lemma 2.1, and taking into account that $\lambda\neq 0,\
\lambda + 2\lambda't\neq 0$, we obtain that $\alpha_1,\ \alpha_2,\
\alpha_3,\ \beta_1,\ \beta_2, \ \beta_3$ must vanish.

\vskip2mm

Solving the equations $\alpha_1=0$, $\alpha_2=0,$ $\alpha_3=0$
with respect to $\rho$ we get the same value of $\rho$, which is
quite long and we shall not write here.

\vskip2mm

Next, from $\beta_1=0$, $\beta_2=0,$ and $\beta_3=0$, we obtain
another three values for $\rho,$ which we denote respectively by
$\rho_1,$ $\rho_2,$ and $\rho_3$. This values must coincide with
$\rho$.

\vskip2mm
When we impose the conditions $\rho_2-\rho=0,$
$\rho_3-\rho=0,$ we obtain two equations:
\begin{eqnarray}\label{ec1}
\begin{array}{c}
  (a_1^2 + a_1^2a_3^2 - 4a_1a_1't - 4a_1a_1'a_3^2t + 4a_1^2a_3a_3't
+ 4a_1'^2t^2 + 4a_1'^2a_3^2t^2- \\
 - 8a_1a_1'a_3a_3't^2 + 4a_1^2a_3'^2t^2)(A n +B)/N_1=0 \\
\end{array}
\end{eqnarray}
\begin{equation}\label{ec2}
\begin{array}{c}
 (a_1^3a_3 - 2a_1^2a_1'a_3t + 2a_1^3a_3't + 2a_1a_3ct +
2a_1a_3^3ct-4a_1'a_3ct^2 - 4a_1'a_3^3ct^2- \\
 - 4a_1a_3'ct^2 + 4a_1a_3^2a_3'ct^2)(A n +B)/N_2 =0\\
\end{array}
\end{equation}
where the expressions of $A,\ B,\ N_1,\ N_2$ are quite long,
depending on $a_1,\ a_3,\ \lambda$, and their derivatives.

\vskip2mm
Let us study the first parenthesis from (\ref{ec1}) and
(\ref{ec2}), namely
$$
E=a_1^2 + a_1^2a_3^2 - 4a_1a_1't - 4a_1a_1'a_3^2t +
4a_1^2a_3a_3't+
$$
$$
+4a_1'^2t^2 + 4a_1'^2a_3^2t^2- 8a_1a_1'a_3a_3't^2 +
4a_1^2a_3'^2t^2,
$$

$$
F=a_1^3a_3 - 2a_1^2a_1'a_3t + 2a_1^3a_3't + 2a_1a_3ct +
2a_1a_3^3ct-
$$
$$
-4a_1'a_3ct^2 - 4a_1'a_3^3ct^2- 4a_1a_3'ct^2 + 4a_1a_3^2a_3'ct^2,
$$

The sign of $E$ may be studied thinking it as a second degree
function of the variable $a_3'$. The associated equation has the
discriminant $\Delta=-(a_1^2t^2(a_1 - 2a_1't)^2)<0, \forall t>0$
and the coefficient of $a_3'^2,$ $4a_1^2t^2>0, \forall t>0$. Thus,
$E>0$ for every $t>0$. If $t=0$, the expression becomes $a_1^2(1 +
a_3^2)>0$. Hence we obtained that $E$ is always positive.

\vskip2mm
Taking into account of the values of $a_3'$ and $a_1'$
from (\ref{integr}) and then multiplying by $\frac{a_1 +
2b_1t}{a_3 + 2b_3t}>0,$ $F=0$ becomes an equation of the second
order with respect to $a_1^2$
\begin{equation}\label{ec3}
~~~~~~~~~~~~~~~~(a_1^2)^2 - 4a_1^2(1 + a_3^2)ct + 4c^2t^2(1 +
a_3^2)^2=0,
\end{equation}
with the discriminant $\Delta=-16a_3^2(1 + a_3^2)c^2t^2<0,\
\forall t>0.$ Thus $F>0, \forall t>0$ and if $t=0,$ $F=a_1^4>0$.

\vskip2mm

Since $E$ and $F$ are always positive, the relations (\ref{ec1}) and
(\ref{ec2}) are fulfilled if and only if $An+B=0$. The obtained
equations does not depend on the dimension $n$ of the base manifold,
so we get that both $A$ and $B$ must vanish.

\vskip2mm

From the condition $A =0$ we get an expression of
$\lambda^{\prime\prime}$, given by
\begin{eqnarray*}
\begin{array}{l}
 \lambda^{\prime \prime} =(a_1^5 (-\lambda^2 (2 a_1^2 a_1^{\prime 2} + a_1^3
a_1^{\prime\prime} + 4 a_1 a_1^\prime a_3^2 c + 8 a_1^2 a_3
a_3^\prime c + 8 a_3^2 c^2 + 8 a_3^4 c^2) -
\\
-2 a_1^2 \lambda \lambda^\prime (a_1 a_1^\prime + 2 a_3^2 c) + 2
a_1^4 \lambda^{\prime 2}) - 2 a_1^4 t (\lambda^2 (-a_1^2
a_1^{\prime 3} - 6 a_1 a_1^{\prime 2} c - 3 a_1^2
a_1^{\prime\prime} c -
\\
-6 a_1 a_1^{\prime 2} a_3^2 c + a_1^2  a_1^{\prime\prime} a_3^2 c
+ 4 a_1^2 a_1^\prime a_3 a_3^\prime c + 2 a_1^3 a_3^{\prime 2} c +
2 a_1^3 a_3  a_3^{\prime\prime} c - 20 a_1^\prime a_3^2 c^2 -
\\
-20 a_1^\prime a_3^4 c^2 - 8 a_1 a_3 a_3^\prime c^2 + 32 a_1 a_3^3
a_3^\prime c^2) + a_1 \lambda \lambda^\prime (-a_1^2 a_1^{\prime
2} + a_1^3  a_1^{\prime\prime} - 6 a_1 a_1^\prime c -
\\
-6 a_1 a_1^\prime a_3^2 c + 12 a_1^2 a_3 a_3^\prime c + 4 a_3^2
c^2 + 4 a_3^4 c^2) +
     2 a_1^3 \lambda^{\prime 2} (a_1 a_1^\prime + 3 c - a_3^2 c)) -
\\
-4 a_1^3 c t^2 (-\lambda^2 (-3 a_1 a_1^{\prime 3} - 3 a_1
a_1^{\prime 3} a_3^2 +
   6 a_1^2 a_1^{\prime 2} a_3 a_3^\prime - 6 a_1^{\prime 2} c - 3 a_1  a_1^{\prime\prime} c -
   \\-20 a_1^{\prime 2} a_3^2 c - 2 a_1  a_1^{\prime\prime} a_3^2 c -
14 a_1^{\prime 2} a_3^4 c +
       a_1  a_1^{\prime\prime} a_3^4 c + 16 a_1 a_1^\prime a_3 a_3^\prime c +
       32 a_1 a_1^\prime a_3^3 a_3^\prime c +
\\
      +4 a_1^2 a_3^{\prime 2} c -  12 a_1^2 a_3^2 a_3^{\prime 2} c + 4 a_1^2 a_3  a_3^{\prime\prime} c -
       4 a_1^2 a_3^3  a_3^{\prime\prime} c) +
     a_1 \lambda \lambda^\prime (3 a_1 a_1^{\prime 2} - 3 a_1^2  a_1^{\prime\prime} +
\\
+3 a_1 a_1^{\prime 2} a_3^2 +
      a_1^2  a_1^{\prime\prime} a_3^2 - 8 a_1^2 a_1^\prime a_3 a_3^\prime + 2 a_1^3 a_3^{\prime 2} +
      2 a_1^3 a_3  a_3^{\prime\prime} + 6 a_1^\prime c - 4 a_1^\prime a_3^2 c -
\\
-10 a_1^\prime a_3^4 c - 16 a_1 a_3 a_3^\prime c + 16 a_1 a_3^3
a_3^\prime c) +
     2 a_1^2 \lambda^{\prime 2} (-3 a_1 a_1^\prime + a_1 a_1^\prime a_3^2 + 2 a_1^2 a_3 a_3^\prime -
     \\
    - 3 c - 2 a_3^2 c + a_3^4 c)) -  8 a_1 c^2 t^3 (-\lambda^2 (3 a_1 a_1^{\prime 3} +
   6 a_1 a_1^{\prime 3} a_3^2 + 3 a_1 a_1^{\prime 3} a_3^4 -
\\
-12 a_1^2 a_1^{\prime 2} a_3 a_3^\prime -
   12 a_1^2 a_1^{\prime 2} a_3^3 a_3^\prime + 12 a_1^3 a_1^\prime a_3^2 a_3^{\prime 2} + 2 a_1^{\prime 2} c +
       a_1  a_1^{\prime\prime} c + 6 a_1^{\prime 2} a_3^2 c +
\\
+ 3 a_1  a_1^{\prime\prime} a_3^2 c +
       6 a_1^{\prime 2} a_3^4 c + 3 a_1  a_1^{\prime\prime} a_3^4 c + 2 a_1^{\prime 2} a_3^6 c +
       a_1  a_1^{\prime\prime} a_3^6 c - 12 a_1 a_1^\prime a_3 a_3^\prime c -
\\
- 24 a_1 a_1^\prime a_3^3 a_3^\prime c -
       12 a_1 a_1^\prime a_3^5 a_3^\prime c - 2 a_1^2 a_3^{\prime 2} c +
       12 a_1^2 a_3^2 a_3^{\prime 2} c + 14 a_1^2 a_3^4 a_3^{\prime 2} c -
\\
- 2 a_1^2 a_3  a_3^{\prime\prime} c -
       4 a_1^2 a_3^3  a_3^{\prime\prime} c - 2 a_1^2 a_3^5  a_3^{\prime\prime} c) +
     a_1 \lambda \lambda^\prime (-3 a_1 a_1^{\prime 2} + 3 a_1^2  a_1^{\prime\prime} + 2 a_1 a_1^{\prime 2} a_3^2 +
\\
+2 a_1^2  a_1^{\prime\prime} a_3^2 + 5 a_1 a_1^{\prime 2} a_3^4 -
a_1^2  a_1^{\prime\prime} a_3^4 +
       8 a_1^2 a_1^\prime a_3 a_3^\prime - 8 a_1^2 a_1^\prime a_3^3 a_3^\prime -
       4 a_1^3 a_3^{\prime 2} -
       \\
       - 4 a_1^3 a_3  a_3^{\prime\prime}+
                            4 a_1^3 a_3^3  a_3^{\prime\prime} - 2 a_1^\prime c - 6 a_1^\prime a_3^2 c -
       6 a_1^\prime a_3^4 c - 2 a_1^\prime a_3^6 c + 4 a_1 a_3 a_3^\prime c +
       \\+8 a_1 a_3^3 a_3^\prime c + 4 a_1 a_3^5 a_3^\prime c) +
        2 a_1^2 \lambda^{\prime 2} (3 a_1 a_1^\prime + 2 a_1 a_1^\prime a_3^2 -
     a_1 a_1^\prime a_3^4 - 4 a_1^2 a_3 a_3^\prime +
\\
+4 a_1^2 a_3^3 a_3^\prime + c + 3 a_3^2 c + 3 a_3^4 c +
       a_3^6 c)) -
       16 c^3 t^4 (\lambda^2 (a_1^\prime + a_1^\prime a_3^2 - 2 a_1 a_3 a_3^\prime)^3
       +\\
      + a_1 \lambda \lambda^\prime (1 + a_3^2) (a_1^{\prime 2} - a_1  a_1^{\prime\prime} + 2 a_1^{\prime 2} a_3^2 -
      2 a_1  a_1^{\prime\prime} a_3^2 +
      a_1^{\prime 2} a_3^4 - a_1  a_1^{\prime\prime} a_3^4 +  2 a_1^2 a_3^{\prime 2} -
      \\
     - 2 a_1^2 a_3^2 a_3^{\prime 2} + 2 a_1^2 a_3  a_3^{\prime\prime} + 2 a_1^2 a_3^3  a_3^{\prime\prime}) +
       2 a_1^2 \lambda^{\prime 2} (1 + a_3^2)^2 (-a_1^\prime - a_1^\prime a_3^2 +
       \\
      + 2 a_1 a_3 a_3^\prime)))/
  (a_1^2 \lambda (a_1^3 - 2 a_1^2 a_1^\prime t - 2 a_1 c t - 2 a_1 a_3^2 c t +
  4 a_1^\prime c t^2 +  4 a_1^\prime a_3^2 c t^2 -
\\
- 8 a_1 a_3 a_3^\prime c t^2)
     (a_1^4 - 4 a_1^2 c t + 4 a_1^2 a_3^2 c t + 4 c^2 t^2 + 8 a_3^2 c^2 t^2 + 4 a_3^4 c^2
     t^2))
\end{array}
\end{eqnarray*}

From $B=0$ we may get the expression of
${\lambda}^{\prime\prime\prime}$, which is quite long and we shall
not present it here.

\vskip2mm
By doing some quite long computations with RICCI, we
prove that the differences $\rho_1-\rho$, $\rho_2-\rho$ and
$\rho_3-\rho$ vanish when we replace the obtained values for
${\lambda}^{\prime\prime}$ and ${\lambda}^{\prime\prime\prime}$.
Hence all the expressions obtained for the constant $\rho$
coincide.

\vskip2mm

Next we have to find the conditions under which the derivative of
${\lambda}^{\prime \prime}$ is equal to ${\lambda}^{\prime \prime
\prime}$:
$$
({\lambda}^{\prime \prime})^{\prime} -{\lambda}^{\prime \prime
\prime} =0.
$$

Computing the above difference, we obtain that its numerator must
vanish:

$$
(a_1^2a_1^{\prime}\lambda + 2a_1c\lambda + 2a_1a_3^2c\lambda +
a_1^3\lambda^\prime -2a_1^{\prime}c{\lambda}t -
$$
$$
-2a_1^{\prime}a_3^2c{\lambda}t + 4a_1a_3a_3^{\prime}c{\lambda}t +
 2a_1c\lambda^{\prime}t +2a_1a_3^2c\lambda^{\prime}t)
$$
$$
(a_1^5a_1^{\prime}\lambda^2 + 2a_1^4a_3^2c\lambda^2 +
a_1^6\lambda\lambda^\prime - a_1^4{a_1^\prime}^2\lambda^2t -
4a_1^3a_1^{\prime}c\lambda^2t - 4a_1^3a_1^{\prime}a_3^2c\lambda^2t
+
$$
$$
+4a_1^4a_3a_3^{\prime}c\lambda^2t -
4a_1^4c\lambda\lambda^{\prime}t +
4a_1^4a_3^2c\lambda\lambda^{\prime}t + a_1^6{\lambda^\prime}^2t +
4a_1^2{a_1^\prime}^2c\lambda^2t^2 +
$$
$$
+4a_1^2{a_1^\prime}^2a_3^2c\lambda^2t^2 -
8a_1^3a_1^{\prime}a_3a_3^{\prime}c\lambda^2t^2 +
4a_1a_1^{\prime}c^2\lambda^2t^2 +
8a_1a_1^{\prime}a_3^2c^2\lambda^2t^2 +
$$
$$
+4a_1a_1^{\prime}a_3^4c^2\lambda^2t^2 -
8a_1^2a_3a_3^{\prime}c^2\lambda^2t^2 -
8a_1^2a_3^3a_3^{\prime}c^2\lambda^2t^2 +
4a_1^2c^2\lambda\lambda^{\prime}t^2 +
$$
$$
+8a_1^2a_3^2c^2\lambda\lambda^{\prime}t^2 +
4a_1^2a_3^4c^2\lambda\lambda^{\prime}t^2 -
4a_1^4c{\lambda^\prime}^2t^2 + 4a_1^4a_3^2c{\lambda^\prime}^2t^2 -
$$
$$
-4{a_1^\prime}^2c^2\lambda^2t^3 -
8{a_1^\prime}^2a_3^2c^2\lambda^2t^3 -
4{a_1^\prime}^2a_3^4c^2\lambda^2t^3 +
16a_1a_1^{\prime}a_3a_3^{\prime}c^2\lambda^2t^3 +
$$
$$
+16a_1a_1^{\prime}a_3^3a_3^{\prime}c^2\lambda^2t^3 -
16a_1^2a_3^2{a_3^\prime}^2c^2\lambda^2t^3 +
4a_1^2c^2{\lambda^\prime}^2t^3 +
8a_1^2a_3^2c^2{\lambda^\prime}^2t^3+
$$
$$
+4a_1^2a_3^4c^2{\lambda^\prime}^2t^3)
$$
$$
(a_1^3 - 2a_1^2a_1^{\prime}t - 2a_1ct - 2a_1a_3^2ct +
4a_1^{\prime}ct^2 + 4a_1^{\prime}a_3^2ct^2 -
8a_1a_3a_3^{\prime}ct^2)= 0.
$$

\vskip2mm

If we replace in the last parenthesis the values of $a_1'$ and
$a_3'$ given by (\ref{integr}) and then we multiply by the
denominator  $(a_1 + 2b_1t)>0,$ we obtain the expression (\ref{ec3})
$$
(a_1^2)^2 - 4a_1^2(1 + a_3^2)ct + 4c^2t^2(1 + a_3^2)^2,
$$
about which we proved that it is always positive.

\vskip2mm
Hence the cases which must be studied are the next two:
$$
I)\ a_1^2a_1^{\prime}\lambda + 2a_1c\lambda + 2a_1a_3^2c\lambda +
a_1^3\lambda^\prime -2a_1^{\prime}c{\lambda}t -
2a_1^{\prime}a_3^2c{\lambda}t + 4a_1a_3a_3^{\prime}c{\lambda}t +
$$
$$
+2a_1c\lambda^{\prime}t +2a_1a_3^2c\lambda^{\prime}t=0,
$$
$$
II)\ a_1^5a_1^{\prime}\lambda^2 + 2a_1^4a_3^2c\lambda^2 +
a_1^6\lambda\lambda^\prime - a_1^4{a_1^\prime}^2\lambda^2t -
4a_1^3a_1^{\prime}c\lambda^2t - 4a_1^3a_1^{\prime}a_3^2c\lambda^2t
+
$$
$$
+4a_1^4a_3a_3^{\prime}c\lambda^2t -
4a_1^4c\lambda\lambda^{\prime}t +
4a_1^4a_3^2c\lambda\lambda^{\prime}t + a_1^6{\lambda^\prime}^2t +
4a_1^2{a_1^\prime}^2c\lambda^2t^2 +
$$
$$
+4a_1^2{a_1^\prime}^2a_3^2c\lambda^2t^2 -
8a_1^3a_1^{\prime}a_3a_3^{\prime}c\lambda^2t^2 +
4a_1a_1^{\prime}c^2\lambda^2t^2 +
8a_1a_1^{\prime}a_3^2c^2\lambda^2t^2 +
$$
$$
+4a_1a_1^{\prime}a_3^4c^2\lambda^2t^2 -
8a_1^2a_3a_3^{\prime}c^2\lambda^2t^2 -
8a_1^2a_3^3a_3^{\prime}c^2\lambda^2t^2 +
4a_1^2c^2\lambda\lambda^{\prime}t^2 +
$$
$$
+8a_1^2a_3^2c^2\lambda\lambda^{\prime}t^2 +
4a_1^2a_3^4c^2\lambda\lambda^{\prime}t^2 -
4a_1^4c{\lambda^\prime}^2t^2 + 4a_1^4a_3^2c{\lambda^\prime}^2t^2 -
$$
$$
-4{a_1^\prime}^2c^2\lambda^2t^3 -
8{a_1^\prime}^2a_3^2c^2\lambda^2t^3 -
4{a_1^\prime}^2a_3^4c^2\lambda^2t^3 +
16a_1a_1^{\prime}a_3a_3^{\prime}c^2\lambda^2t^3 +
$$
$$
+16a_1a_1^{\prime}a_3^3a_3^{\prime}c^2\lambda^2t^3 -
16a_1^2a_3^2{a_3^\prime}^2c^2\lambda^2t^3 +
4a_1^2c^2{\lambda^\prime}^2t^3 +
8a_1^2a_3^2c^2{\lambda^\prime}^2t^3+
$$
$$
+4a_1^2a_3^4c^2{\lambda^\prime}^2t^3=0.
$$

\vskip2mm In the case \emph{I}, we may obtain the following
expression of $\lambda^\prime$
\begin{eqnarray*}
\lambda^\prime=
-\lambda\frac{a_1(a_1a_1^{\prime} + 2c(1+ a_3^2)) -
2ct(a_1^{\prime} + 2a_1^{\prime}a_3^2 -
4a_1a_3a_3^{\prime})}{a_1[a_1^2 + 2ct(1 + a_3^2)]}.
\end{eqnarray*}

Replacing this expression of $\lambda^\prime$ in the first value
obtained for $\rho$, we get a rational function of the density of
energy, t, the coefficients $a_1$ and $a_3$, the proportionality
factor $\lambda$, and the constant sectional curvature, c, of the
base manifold $M$
\begin{eqnarray*}
~~~~~~~~~~~~~~~~~~~~~~~~~~~\rho= \frac{2a_1c(n+1)}{\lambda[a_1^2 +
2ct(1 + a_3^2)]},
\end{eqnarray*}
from which we get the value of $\lambda$
\begin{equation}\label{la}
~~~~~~~~~~~~~~~~~~~~~~~~~~~\lambda = \frac{2a_1c(n+1)}{\rho[a_1^2
+ 2ct(1 + a_3^2)]}.
\end{equation}

Now we may state:

\vskip2mm

\bf{Theorem 3.2.} \it{Let $(M,g)$ be a smooth $n$-dimensional
Riemannian manifold. If $(G,J)$ is a general natural K\" ahlerian
structure on the cotangent bundle $T^*M$ and the parameter
$\lambda$ is expressed by \rm{(\ref{la})}, \it where $\rho$ is a
nonzero real constant, then $(T^*M,G,J)$ is a K\" ahler-Einstein
manifold, i.e. $Ric=\rho G$.}

\vskip2mm {\bf Remark 3.1.} \rm Taking into account of a theorem
from \cite{Druta2}, the expression (\ref{la}) of $\lambda$ implies
that $(T^*M,G,J)$ is a K\" ahlerian manifold of constant
holomorphic sectional curvature $k=\frac{2\rho}{n+1}$. \vskip2mm

\vskip2mm {\bf Example 3.1.} The K\"ahler-Einstein structure on
$T^*M$, from the paper \cite{OprPor1} by  Oproiu and Poro\c sniuc,
may be obtained from the theorem 3.2, as a particular case.  If in
the expression (\ref{la}) we impose the condition $a_3=0,$ we get
the same expression of $\lambda$ obtained in \cite{OprPor1}, in
the case of the natural structure of diagonal lifted type on the
cotangent bundle $T^*M$ of a Riemannian manifold $(M,g)$.

\vskip3mm

\rm In the case \emph{II}, we obtain
$$
a_1^2t(a_1^4 - 4a_1^2ct + 4a_1^2a_3^2ct + 4c^2t^2 + 8a_3^2c^2t^2 +
4a_3^4c^2t^2){\lambda^\prime}^2 + a_1^2(a_1^4 - 4a_1^2ct +
$$
$$
+4a_1^2a_3^2ct + 4c^2t^2 + 8a_3^2c^2t^2 +
4a_3^4c^2t^2)\lambda^{\prime}\lambda + (a_1^5a_1^\prime +
2a_1^4a_3^2c - a_1^4{a_1^\prime}^2t - 4a_1^3a_1^{\prime}ct -
$$
$$
-4a_1^3a_1^{\prime}a_3^2ct + 4a_1^4a_3a_3^{\prime}ct +
4a_1^2{a_1^\prime}^2ct^2 +
  4a_1^2{a_1^\prime}^2a_3^2ct^2 - 8a_1^3a_1^{\prime}a_3a_3^{\prime}ct^2 +
$$
$$
 + 4a_1a_1^{\prime}c^2t^2 + 8a_1a_1^{\prime}a_3^2c^2t^2 + 4a_1a_1^{\prime}a_3^4c^2t^2 -
  8a_1^2a_3a_3^{\prime}c^2t^2 - 8a_1^2a_3^3a_3^{\prime}c^2t^2 -
$$
$$
 - 4{a_1^\prime}^2c^2t^3 - 8{a_1^\prime}^2a_3^2c^2t^3 - 4{a_1^\prime}^2a_3^4c^2t^3 +
  16a_1a_1^{\prime}a_3a_3^{\prime}c^2t^3 + 16a_1a_1^{\prime}a_3^3a_3^{\prime}c^2t^3 -
$$
$$
  - 16a_1^2a_3^2{a_3^\prime}^2c^2t^3)\lambda^2 =0.
$$

This is a homogeneous equation of second order in $\lambda^\prime$
and $\lambda $ and it may be solved with respect to
$\frac{\lambda^\prime}{\lambda}$. Then we obtain two expressions
for $\lambda^\prime$
$$
\lambda^\prime =\lambda (\pm\frac{1}{2 t} + \frac{a_1^3 - 2 a_1^2
a_1^\prime t - 2 a_1 c t - 2 a_1 a_3^2 c t + 4 a_1^\prime c t^2 +
4 a_1^\prime a_3^2 c t^2 - 8 a_1 a_3 a_3^\prime c t^2}{
   {2 a_1 t \sqrt{a_1^4 - 4 a_1^2 c t + 4 a_1^2 a_3^2 c t + 4 c^2 t^2 + 8 a_3^2 c^2 t^2 + 4 a_3^4 c^2
   t^2}}}).
$$

When we replace this expression of ${\lambda}^\prime$ and its
derivative ${\lambda}^{\prime\prime}$ in the first value of
$\rho$, we obtain
$$
~~~~~~~\rho=\frac{n (a_1^2 + 2 c t + 2 a_3^2 c t \pm
            \sqrt{a_1^4 - 4 a_1^2 c t + 4 a_1^2 a_3^2 c t + 4 c^2 t^2 +
8 a_3^2 c^2 t^2 + 4 a_3^4 c^2 t^2})}{4 a_1 \lambda t}.
$$

In this case $\rho$ is defined on the set $T_0M\subset TM$ of the
nonzero cotangent vectors to $M$, and the value of $\lambda$ is
given by
\begin{equation}\label{la1}
\lambda=\frac{n (a_1^2 + 2 c t + 2 a_3^2 c t \pm \sqrt{a_1^4 - 4
a_1^2 c t + 4 a_1^2 a_3^2 c t + 4 c^2 t^2 + 8 a_3^2 c^2 t^2 + 4
a_3^4 c^2 t^2})}{4 a_1 \rho t}.
\end{equation}

Now we may formulate the next theorem:

\vskip2mm

\bf{Theorem 3.3.} \it{Let $(G,J)$ be a general natural K\"
ahlerian structure on the cotangent bundle $T^*M$ of a smooth
$n$-dimensional Riemannian manifold. If the parameter $\lambda$ is
expressed by \rm{(\ref{la1})}, \it where $\rho$ is a nonzero real
constant, then $(G,J)$ is a K\"ahler-Einstein structure on the
bundle $T_0^*M$, of nonzero cotangent vectors to $M$, i.e.
$Ric=\rho G$.}

\vskip2mm \bf{Aknowledgement.}\rm The author expresses her
gratitude to Professor Oproiu, her PhD advisor, for the
mathematical conversations throughout this work, for the
scientific support and the techniques learned being his PhD
student.

This work was partially supported by the Grant TD-158/2007,
CNCSIS, Ministerul  Educa\c tiei \c si Cercet\u arii, Rom\^ania.

\vskip1.5mm

\begin{minipage}{2.5in}
\begin{flushleft}
\it Author's adress:\vskip2mm

\rm Simona-Luiza Dru\c t\u a\\
\rm Faculty of Mathematics, \\
University "Al.I. Cuza", Ia\c si, \\
RO-700 506, Romania.\\
e-mail: simonadruta@yahoo.com
\end{flushleft}
\end{minipage}

\end{document}